\documentclass{article}

\usepackage[top=20truemm,bottom=20truemm,left=20truemm,right=20truemm]{geometry}
\usepackage{authblk}
\usepackage{abstract}
\usepackage{amsmath}
\usepackage{amssymb}
\usepackage{hyperref}
\usepackage[dvipdfmx]{graphicx}
\usepackage{color}
\usepackage{diagbox}

\author[1]{Takuya Tsuchiya\thanks{\href{mailto:tatsuchi@eco.meijigakuin.ac.jp}{\nolinkurl{tatsuchi@eco.meijigakuin.ac.jp}}}}
\affil[1]{%
  Faculty of Economics, Meiji Gakuin University, Japan
}

\title{Structure-preserving numerical calculation of wave equation for a vector
  field}
\begin{document}

\maketitle

\begin{abstract}
  For the Proca equation, which is a wave equation for a vector field, we derive
  the canonical formulation including constraints from the Stueckelberg action
  and propose discrete equations with a structure-preserving scheme for
  conserving the constraints at the discrete level.
  Numerical simulations are performed using these discrete equations and other
  discrete equations with a standard scheme.
  We show the results obtained using the structure-preserving scheme and provide
  more accurate and stable numerical solutions.
\end{abstract}

\section{Introduction}

The wave equation is one of the dynamical equations that describe various
natural phenomena, including the propagation of electromagnetic and
gravitational waves.
These equations often have constraints.
When performing numerical calculations using these equations with constraints,
the constraints often become unsatisfiable for time evolution owing to the
accumulation of numerical errors.
Thus, to carry out high-precision and stable numerical solutions, it is
necessary to perform numerical calculations using discrete equations with a
structure that does not accumulate numerical errors.

We have been studying high-precision and stable numerical calculations for a
scalar field \cite{Tsuchiya-Nakamura-2019-JCAM, Tsuchiya-Nakamura-2022-ISAAC,
  Tsuchiya-Nakamura-2023-JSIAML, Tsuchiya-Nakamura-2025}, which is one of the
wave equations.
In this paper, we investigate the Proca equation \cite{Proca}, which is a
vector field equation known as the Maxwell equation with a mass term.

Indices such as $(i, j, \dots)$ and $(\mu,\nu,\dots)$ run from $1$ to $n$ and
$0$ to $n$, respectively, where $n$ is the spatial dimension.
We use the Einstein convention of summation of repeated up--down indices in this
paper.

\section{Stueckelberg action and canonical formulation}

The Proca equation is given as
\begin{align}
  \partial_\nu\partial^\nu A^\mu - \partial^\mu \partial_\nu A^\nu
  + \dfrac{c^2m^2}{\hbar^2}A^\mu=0,
  \label{eq:ProcaEq}
\end{align}
where $A^\mu$ is the dynamical variable, $m$ is the mass, $c$ is the speed of
light, and $\hbar$ is the Dirac constant.
The divergence of \eqref{eq:ProcaEq} gives the equation
\begin{align}
  \partial_\mu A^\mu =0.
  \label{eq:ProcaEq2}
\end{align}

We give the Lagrangian density, which is called the Stueckelberg action \cite{%
  Stueckelberg, StueckelbergAction}, as follows:
\begin{align}
  \mathcal{L}
  &= \dfrac{1}{4p_1}(\partial_\mu A_\nu -\partial_\nu A_\mu)
  (\partial^\mu A^\nu -\partial^\nu A^\mu)
  + \dfrac{p_2}{2}A^\mu A_\mu
  + \dfrac{1}{2\lambda}(\partial_\mu A^\mu)(\partial_\nu A^\nu),
  \label{eq:PSform}
\end{align}
where $p_1\neq 0$ and $p_2\in\mathbb{R}$ are constant values, and
$\lambda\neq 0$.
Using the Euler--Lagrange equation of \eqref{eq:PSform}, we can derive the Proca
equations \eqref{eq:ProcaEq} and \eqref{eq:ProcaEq2} if we set $p_1=1$ and
$p_2=c^2m^2/\hbar^2$.
The Hamiltonian density of \eqref{eq:PSform} is given by
\begin{align}
  \mathcal{H}
  &= \dfrac{\lambda}{2}\Pi_0\Pi_0
  - \Pi_0(\partial_mA^m)
  + \dfrac{p_1}{2}\Pi_m\Pi^m
  - \Pi^m(\partial_mA^0)
  + \dfrac{1}{2p_1}(\partial_mA_n)(\partial^mA^n)
  - \dfrac{1}{2p_1}(\partial_nA^m)(\partial_mA^n)
  + \dfrac{p_2}{2} A^0A^0
  \nonumber\\
  &\quad
  - \dfrac{p_2}{2} A_m A^m,
  \label{eq:ProcaHamiltonianDensity}
\end{align}
where $\Pi_0$ and $\Pi_i$ are the canonical momenta of $A^0$ and $A^i$,
respectively.
Then, the canonical equations of \eqref{eq:ProcaHamiltonianDensity} are
\begin{align}
  \partial_0A^0
  &= \lambda\Pi_0
  - \partial_iA^i,
  \label{eq:PSA0Evo1}\\
  \partial_0\Pi_0
  &= -p_2A^0
  -\partial_i\Pi^i,
  \label{eq:PSPi0Evo1}\\
  \partial_0A^i
  &= p_1\Pi^i
  - \partial^iA^0,
  \label{eq:PSAiEvo1}\\
  \partial_0\Pi_i
  &=
  p_2A_i
  - \partial_i\Pi_0
  + \dfrac{1}{p_1}(\partial_j\partial^jA_i-\partial_i\partial_jA^j),
  \label{eq:PSPiiEvo1}
\end{align}
where $\partial_0= (1/c)\partial_t$.
The constraint equation is derived from the variation of \eqref{%
  eq:ProcaHamiltonianDensity} with respect to $\lambda$, such as
\begin{align}
  \mathcal{C}_1:= \Pi_0\approx0,
  \label{eq:Const1}
\end{align}
since  $\lambda$ is the gauge variable, where $\approx$ is called the weak
equality, which means equality analytically, but is used when the equality is
not satisfied owing to numerical errors during the numerical calculation.
Note that $\Pi_0$ is both a dynamical variable and a constraint variable.
In addition, since the time derivative of \eqref{eq:Const1} should be also the
constraint equation, the following condition is required:
\begin{align}
  \mathcal{C}_2
  &:= \partial_0\mathcal{C}_1
  = -p_2A^0 - \partial_i\Pi^i
  \approx 0.
  \label{eq:Const2}
\end{align}
This equation is equivalent to Gauss's law in the Maxwell equation.
The time derivative of $\mathcal{C}_2$ is calculated as
\begin{align}
  \partial_0\mathcal{C}_{2}
  &= -p_2\lambda \mathcal{C}_{1}
  + \partial_i\partial^i \mathcal{C}_{1}.
  \label{eq:timederiveConst2}
\end{align}
Thus, $\mathcal{C}_1$ and $\mathcal{C}_2$ are conserved in time evolutions if
$\mathcal{C}_1\approx 0$ and $\mathcal{C}_2\approx 0$ are satisfied at the
initial time.

If we perform the Fourier analysis of \eqref{eq:Const2} and \eqref{%
  eq:timederiveConst2},
the constraints are expressed as $\mathcal{C}_{i}:=
\int \hat{\mathcal{C}}_{i}(t,h^j)\exp({\rm i} h_m x^m)d^nk$,
\begin{align}
  \left\{
  \begin{array}{l}
    \partial_0\hat{\mathcal{C}}_{1}
    = \hat{\mathcal{C}}_{2},\\
    \partial_0\hat{\mathcal{C}}_{2}
    = -p_2\lambda \hat{\mathcal{C}}_{1}
    - h_ih^i\hat{\mathcal{C}}_{1},\\
  \end{array}
  \right.
  \Leftrightarrow\quad
  \partial_0\begin{pmatrix}
  \hat{\mathcal{C}}_{1}\\
  \hat{\mathcal{C}}_{2}
  \end{pmatrix}
  =
  \begin{pmatrix}
    0 & 1\\
    -p_2\lambda -h^ih_i & 0
  \end{pmatrix}
  \begin{pmatrix}
  \hat{\mathcal{C}}_{1}\\
  \hat{\mathcal{C}}_{2}
  \end{pmatrix},
  \label{eq:CP1}
\end{align}
where $h^i$ is the wave vector.
The eigenvalues of the coefficient matrix of the right-hand side of \eqref{%
  eq:CP1} are $\pm{\rm i}\sqrt{p_2\lambda+h_ih^i}$.
If $p_2\lambda+h_ih^i>0$, there are no modes that increase the magnitudes of the
constraints since both of the eigenvalues are pure imaginary numbers.
On the other hand, purely imaginary eigenvalues do not reduce the numerical
errors that accumulate on the constraint errors during numerical calculations,
and once it becomes unstable, the numerical simulations remain unstable.

The time derivative of the total Hamiltonian
$H_C:=\int_{\mathbb{R}^n}\mathcal{H}d^nx$ is calculated as
\begin{align}
  \partial_0H_C
  &=
  \int_{\mathbb{R}^n}
  \dfrac{1}{2}(\partial_0\lambda)(\mathcal{C}_{1})^2d^nx
  + [\mathrm{Boundary\,\,Terms}].
\end{align}
Thus, $H_C$ is conserved if $\mathcal{C}_1\approx0$ is satisfied and appropriate
boundary conditions are imposed.

\section{Structure-preserving discrete equations}

The discrete equations of \eqref{eq:ProcaHamiltonianDensity}--\eqref{eq:Const2}
are
\begin{align}
  \mathcal{H}^{(\ell)}_{(\boldsymbol{k})}
  &:=
  \dfrac{1}{2}\lambda^{(\ell)}_{(\boldsymbol{k})}
  \Pi_0{}^{(\ell)}_{(\boldsymbol{k})}\Pi_0{}^{(\ell)}_{(\boldsymbol{k})}
  - \Pi_0{}^{(\ell)}_{(\boldsymbol{k})}
  (\widehat{\delta}^{\langle1\rangle}_mA^m{}^{(\ell)}_{(\boldsymbol{k})})
  + \dfrac{p_1}{2}\Pi^m{}^{(\ell)}_{(\boldsymbol{k})}
  \Pi_m{}^{(\ell)}_{(\boldsymbol{k})}
  - \Pi^m{}^{(\ell)}_{(\boldsymbol{k})}
  (\widehat{\delta}^{\langle1\rangle}_mA^0{}^{(\ell)}_{(\boldsymbol{k})})
  + \dfrac{1}{2p_1}
  (\widehat{\delta}^{\langle1\rangle}_mA^n{}^{(\ell)}_{(\boldsymbol{k})})
  (\widehat{\delta}^{\langle1\rangle}{}^mA_n{}^{(\ell)}_{(\boldsymbol{k})})
  \nonumber\\
  &\quad
  + \dfrac{p_2}{2} A^0{}^{(\ell)}_{(\boldsymbol{k})}
  A^0{}^{(\ell)}_{(\boldsymbol{k})}
  - \dfrac{1}{2p_1}
  (\widehat{\delta}^{\langle1\rangle}_nA^m{}^{(\ell)}_{(\boldsymbol{k})})
  (\widehat{\delta}^{\langle1\rangle}_mA^n{}^{(\ell)}_{(\boldsymbol{k})})
  - \dfrac{p_2}{2}A^m{}^{(\ell)}_{(\boldsymbol{k})}
  A_m{}^{(\ell)}_{(\boldsymbol{k})},
  \label{eq:HamitonianDens_SPS}
\end{align}
\begin{align}
  \dfrac{A^{0}{}^{(\ell+1)}_{(\boldsymbol{k})}
    -A^{0}{}^{(\ell)}_{(\boldsymbol{k})}}{c\Delta t}
  &:= \dfrac{1}{4}(\lambda{}^{(\ell+1)}_{(\boldsymbol{k})}
  + \lambda{}^{(\ell)}_{(\boldsymbol{k})})
  (\Pi_{0}{}^{(\ell+1)}_{(\boldsymbol{k})}
  + \Pi_{0}{}^{(\ell)}_{(\boldsymbol{k})})
  - \dfrac{1}{2}
  \widehat{\delta}^{\langle1\rangle}_i
  (A^i{}^{(\ell+1)}_{(\boldsymbol{k})}+A^i{}^{(\ell)}_{(\boldsymbol{k})}),
  \label{eq:A0evolve_SPS}\\
  \dfrac{\Pi_{0}{}^{(\ell+1)}_{(\boldsymbol{k})}
    -\Pi_{0}{}^{(\ell)}_{(\boldsymbol{k})}}{c\Delta t}
  &:= - \dfrac{p_2}{2}(A^{0}{}^{(\ell+1)}_{(\boldsymbol{k})}
  + A^{0}{}^{(\ell)}_{(\boldsymbol{k})})
  - \dfrac{1}{2}\widehat{\delta}^{\langle1\rangle}_i
  (\Pi^i{}^{(\ell+1)}_{(\boldsymbol{k})}
  + \Pi^i{}^{(\ell)}_{(\boldsymbol{k})}),
  \label{eq:Pi0evolve_SPS}\\
  \dfrac{A^{i}{}^{(\ell+1)}_{(\boldsymbol{k})}
    -A^{i}{}^{(\ell)}_{(\boldsymbol{k})}}{c\Delta t}
  &:= \dfrac{p_1}{2}(\Pi^i{}^{(\ell+1)}_{(\boldsymbol{k})}
  + \Pi^i{}^{(\ell)}_{(\boldsymbol{k})})
  - \dfrac{1}{2}\widehat{\delta}^{\langle1\rangle}{}^i
  (A^0{}^{(\ell+1)}_{(\boldsymbol{k})}+A^0{}^{(\ell)}_{(\boldsymbol{k})}),
  \label{eq:Aievolve_SPS}\\
  \dfrac{\Pi_{i}{}^{(\ell+1)}_{(\boldsymbol{k})}
    -\Pi_{i}{}^{(\ell)}_{(\boldsymbol{k})}}{c\Delta t}
  &:=
  - \dfrac{1}{2}\widehat{\delta}^{\langle1\rangle}_{i}
  (\Pi_0{}^{(\ell+1)}_{(\boldsymbol{k})} + \Pi_0{}^{(\ell)}_{(\boldsymbol{k})})
  + \dfrac{p_2}{2}(A_{i}{}^{(\ell+1)}_{(\boldsymbol{k})}
  + A_{i}{}^{(\ell)}_{(\boldsymbol{k})})
  + \dfrac{\widehat{\delta}^{\langle1\rangle}_m
  \widehat{\delta}^{\langle1\rangle}{}^m(A_i{}^{(\ell+1)}_{(\boldsymbol{k})}
  + A_i{}^{(\ell)}_{(\boldsymbol{k})})}{2p_1}
  \nonumber\\
  &\quad
  - \dfrac{\widehat{\delta}^{\langle1\rangle}_i
  \widehat{\delta}^{\langle1\rangle}_m(A^m{}^{(\ell+1)}_{(\boldsymbol{k})}
  + A^m{}^{(\ell)}_{(\boldsymbol{k})})}{2p_1},
  \label{eq:Piievolve_SPS}\\
  \mathcal{C}_{1}{}^{(\ell)}_{(\boldsymbol{k})}
  &:= \Pi_{0}{}^{(\ell)}_{(\boldsymbol{k})},
  \label{eq:ConstC1_SPS}
  \\
  \mathcal{C}_{2}{}^{(\ell)}_{(\boldsymbol{k})}
  &:= -p_2A^{0}{}^{(\ell)}_{(\boldsymbol{k})}
  - \widehat{\delta}^{\langle1\rangle}_i\Pi^i{}^{(\ell)}_{(\boldsymbol{k})},
  \label{eq:ConstC2_SPS}
\end{align}
respectively, where ${}^{(\ell)}$ means the time index, ${}_{(\boldsymbol{k})}$
means the space index, and $\boldsymbol{k}=(k_1,\dots,k_n)$.
$\widehat{\delta}^{\langle1\rangle}_i$ is the well-known first-order central
difference operator for the $x^i$ dimension (cf. \cite{%
  Tsuchiya-Nakamura-2023-JSIAML}).
We call the discrete equations \eqref{eq:HamitonianDens_SPS}--\eqref{%
  eq:ConstC2_SPS} as the structure-preserving scheme (SPS) system.

The time derivatives of \eqref{eq:ConstC1_SPS} and \eqref{eq:ConstC2_SPS} are
calculated as
\begin{align}
  &\dfrac{\mathcal{C}_{1}{}^{(\ell+1)}_{(\boldsymbol{k})}
    - \mathcal{C}_{1}{}^{(\ell)}_{(\boldsymbol{k})}}{c\Delta t}
  =\dfrac{1}{2}(\mathcal{C}_{2}{}^{(\ell+1)}_{(\boldsymbol{k})}
  + \mathcal{C}_{2}{}^{(\ell)}_{(\boldsymbol{k})}),
  \label{eq:DiscreteDeriveC1}
  \\
  &\dfrac{\mathcal{C}_{2}{}^{(\ell+1)}_{(\boldsymbol{k})}
    - \mathcal{C}_{2}{}^{(\ell)}_{(\boldsymbol{k})}}{c\Delta t}
  =
  - \dfrac{p_2}{4}(\lambda{}^{(\ell+1)}_{(\boldsymbol{k})}
  + \lambda{}^{(\ell)}_{(\boldsymbol{k})})
  (\mathcal{C}_{1}{}^{(\ell+1)}_{(\boldsymbol{k})}
  + \mathcal{C}_{1}{}^{(\ell)}_{(\boldsymbol{k})})
  + \dfrac{1}{2}\widehat{\delta}^{\langle1\rangle}_{i}
  \widehat{\delta}^{\langle1\rangle}{}^i
  (\mathcal{C}_{1}{}^{(\ell+1)}_{(\boldsymbol{k})}
  + \mathcal{C}_{1}{}^{(\ell)}_{(\boldsymbol{k})}),
  \label{eq:DiscreteDeriveC2}
\end{align}
respectively.
From these results, $\mathcal{C}_1{}^{(\ell+1)}_{(\boldsymbol{k})}\approx0$
and $\mathcal{C}_2{}^{(\ell+1)}_{(\boldsymbol{k})}\approx0$ are satisfied if
$\mathcal{C}_1{}^{(\ell)}_{(\boldsymbol{k})}\approx0$
and $\mathcal{C}_2{}^{(\ell)}_{(\boldsymbol{k})}\approx0$ are satisfied.

The discrete total Hamiltonian is defined as
\begin{align}
  {H}_C{}^{(\ell)}
  := \sum_{D}\mathcal{H}{}^{(\ell)}_{(\boldsymbol{k})}\Delta V,
  \label{eq:discreteTotalH}
\end{align}
where $D=\{3\le k_1\le N_1+2,\dots,3\le k_n\le N_n+2\}$, $N_i$ is the number of
grids of $x^i$, and $\Delta V:=\prod_{1\leq i\leq n}\Delta x^i$.
Then, the discrete derivative of $H_C$ in time is
\begin{align}
  &\frac{H_C{}^{(\ell+1)} - H_C{}^{(\ell)}}{c\Delta t}
  = \sum_{D}\dfrac{\mathcal{H}{}^{(\ell+1)}_{(\boldsymbol{k})}
    - \mathcal{H}{}^{(\ell)}_{(\boldsymbol{k})}}{c\Delta t}\Delta V
  = \sum_{D}\dfrac{\Delta V}{4}\dfrac{\lambda^{(\ell+1)}_{(\boldsymbol{k})}
    -\lambda^{(\ell)}_{(\boldsymbol{k})}}{c\Delta t}
  \{(\mathcal{C}_{1}{}^{(\ell+1)}_{(\boldsymbol{k})})^2
  + (\mathcal{C}_{1}{}^{(\ell)}_{(\boldsymbol{k})})^2\}
  + [\text{Boundary\,\,Terms}].
  \label{eq:DiscTotalHtimeDerive}
\end{align}
If we set appropriate boundary conditions, and $\mathcal{C}_1{}^{(\ell)}_{(\boldsymbol{k})}\approx 0$ and $\mathcal{C}_2{}^{(\ell)}_{(\boldsymbol{k})}\approx0$,
then $H_C{}^{(\ell)}$ is preserved in time evolutions.

\section{Numerical tests}

In this section, we perform some simulations to confirm that \eqref{%
  eq:DiscreteDeriveC1}--\eqref{eq:DiscTotalHtimeDerive} are satisfied
numerically.
The initial conditions are set as
\begin{align}
  A^0&=-\dfrac{2a\pi\{\cos (2\pi (x+y))+\sin (2\pi (x+y))\}}{%
    \sqrt{8\pi^2-p_1p_2}},\\
  \Pi_0&=0,\\
  A^1 &= a\cos(2\pi (x+y)),\\
  A^2 &= a\sin(2\pi (x+y)),\\
  A^3 &= a\cos(4\pi (x+y)),\\
  \Pi_1&= \dfrac{-4a\pi^2\cos(2\pi(x+y))}{p_1\sqrt{8\pi^2-p_1p_2}}
  + \dfrac{a(p_1p_2-4\pi^2)\sin(2\pi (x+y))}{p_1\sqrt{8\pi^2-p_1p_2}}
  ,
\end{align}
\begin{align}
  \Pi_2 &= \dfrac{4a\pi^2\sin(2\pi(x+y))}{p_1\sqrt{8\pi^2-p_1p_2}}
  + \dfrac{a(-p_1p_2+4\pi^2)\cos(2\pi (x+y))}{p_1\sqrt{8\pi^2-p_1p_2}},\\
  \Pi_3 &=-\dfrac{2a\sqrt{8\pi^2-p_1p_2}}{p_1}\sin(4\pi (x+y)),\\
  \lambda&=0.01,
\end{align}
where $a=1$ and $-1/2\leq x,y\leq 1/2$.
The boundary is periodic.
The parameters are $c=p_1=p_2=1$.
The spatial dimension is $n=3$.
The spatial grid ranges are $\Delta x=\Delta y=1/50, 1/100, 1/200$.
The CFL condition is $1/4$.
The simulation time is $0\leq t\leq 80$.

To compare the numerical results with the SPS results, we propose other discrete
equations, namely,
\begin{align}
  \mathcal{H}^{(\ell)}_{(\boldsymbol{k})}
  &:=
  \dfrac{1}{2}\lambda^{(\ell)}_{(\boldsymbol{k})}
  \Pi_0{}^{(\ell)}_{(\boldsymbol{k})}\Pi_0{}^{(\ell)}_{(\boldsymbol{k})}
  - \Pi_0{}^{(\ell)}_{(\boldsymbol{k})}
  (\widehat{\delta}^{\langle1\rangle}_mA^m{}^{(\ell)}_{(\boldsymbol{k})})
  + \dfrac{p_1}{2}
  \Pi_m{}^{(\ell)}_{(\boldsymbol{k})}\Pi^m{}^{(\ell)}_{(\boldsymbol{k})}
  + \dfrac{1}{4p_1}
  (\widehat{\delta}^{+}_mA^n{}^{(\ell)}_{(\boldsymbol{k})})
  (\widehat{\delta}^{+}{}^mA_n{}^{(\ell)}_{(\boldsymbol{k})})
  \nonumber\\
  &\quad
  + \dfrac{1}{4p_1}
    (\widehat{\delta}^{-}_mA^n{}^{(\ell)}_{(\boldsymbol{k})})
    (\widehat{\delta}^{-}{}^mA_n{}^{(\ell)}_{(\boldsymbol{k})})
  + \dfrac{p_2}{2} A^0{}^{(\ell)}_{(\boldsymbol{k})}
  A^0{}^{(\ell)}_{(\boldsymbol{k})}
  - \dfrac{1}{4p_1}
  (\widehat{\delta}^{+}_nA^m{}^{(\ell)}_{(\boldsymbol{k})})
  (\widehat{\delta}^{+}_mA^n{}^{(\ell)}_{(\boldsymbol{k})})
  - \dfrac{p_2}{2}A^m{}^{(\ell)}_{(\boldsymbol{k})}
  A_m{}^{(\ell)}_{(\boldsymbol{k})}
  \nonumber\\
  &\quad
  - \dfrac{1}{4p_1}
  (\widehat{\delta}^{-}_nA^m{}^{(\ell)}_{(\boldsymbol{k})})
  (\widehat{\delta}^{-}_mA^n{}^{(\ell)}_{(\boldsymbol{k})})
  - \Pi^m{}^{(\ell)}_{(\boldsymbol{k})}
  (\widehat{\delta}^{\langle1\rangle}_mA^0{}^{(\ell)}_{(\boldsymbol{k})}),
  \label{eq:HamiltonianDens_CS}
\end{align}
\begin{align}
  \dfrac{\Pi_{i}{}^{(\ell+1)}_{(\boldsymbol{k})}
    - \Pi_{i}{}^{(\ell)}_{(\boldsymbol{k})}}{c\Delta t}
  &:= \dfrac{p_2}{2}(A_{i}{}^{(\ell+1)}_{(\boldsymbol{k})}
  + A_{i}{}^{(\ell)}_{(\boldsymbol{k})})
  - \dfrac{1}{2}\widehat{\delta}^{\langle1\rangle}_{i}
  (\Pi_0{}^{(\ell+1)}_{(\boldsymbol{k})}+\Pi_0{}^{(\ell)}_{(\boldsymbol{k})})
  \nonumber\\
  &\quad
  + \dfrac{\{\widehat{\delta}^{\langle2\rangle}{}^m{}_{m}
  (A_i{}^{(\ell+1)}_{(\boldsymbol{k})} + A_i{}^{(\ell)}_{(\boldsymbol{k})})
  - \widehat{\delta}^{\langle2\rangle}_{im}
  (A^m{}^{(\ell+1)}_{(\boldsymbol{k})} + A^m{}^{(\ell)}_{(\boldsymbol{k})})\}}{%
    2p_1},
  \label{eq:Piievolve_SS}
\end{align}
\eqref{eq:A0evolve_SPS}--\eqref{eq:Aievolve_SPS}, and \eqref{eq:ConstC1_SPS}--%
\eqref{eq:ConstC2_SPS}.
The symbols $\widehat{\delta}^{+}_{i}$ and $\widehat{\delta}^{-}_{i}$ are the
well-known space forward and backward difference operators, respectively.
The symbol $\widehat{\delta}^{\langle2\rangle}_{ij}$ is the second-order
central difference operator for $x^i$ and $x^j$ dimensions, which is also
well-known if $i=j$ (cf. \cite{Tsuchiya-Nakamura-2022-ISAAC}).
We call these discrete equations the standard scheme (SS) system since the
discrete operator for the second-order derivative in \eqref{eq:Piievolve_SS} is
the standard expression.
For SS, the time derivative of $\mathcal{C}_2{}^{(\ell)}_{(\boldsymbol{k})}$ is
\begin{align}
  \dfrac{\mathcal{C}_2{}^{(\ell+1)}_{(\boldsymbol{k})}
    - \mathcal{C}_2{}^{(\ell)}_{(\boldsymbol{k})}}{c\Delta t}
  &= \text{r.h.s.\,\,of\,\,\eqref{eq:DiscreteDeriveC2}}
  - \dfrac{(\widehat{\delta}^{\langle1\rangle}{}_i
  \widehat{\delta}^{\langle2\rangle}{}^m{}_m
  - \widehat{\delta}^{\langle1\rangle}{}_m
  \widehat{\delta}^{\langle2\rangle}{}^m{}_i)
  (A^i{}^{(\ell+1)}_{(\boldsymbol{k})}+A^i{}^{(\ell)}_{(\boldsymbol{k})})%
  }{2p_1},
  \label{eq:Const2Derive_SS}
\end{align}
and $\widehat{\delta}^{\langle1\rangle}{}_i
\widehat{\delta}^{\langle2\rangle}{}^m{}_m
- \widehat{\delta}^{\langle1\rangle}{}_m
\widehat{\delta}^{\langle2\rangle}{}^m{}_i\neq 0$ if $i\neq m$.
Thus, generally, $\mathcal{C}_2{}^{(\ell)}_{(\boldsymbol{k})}$ is not satisfied
with SS.

\begin{figure}[htbp]
  \centering
  \begin{minipage}{0.32\hsize}
    \includegraphics[width=\hsize]{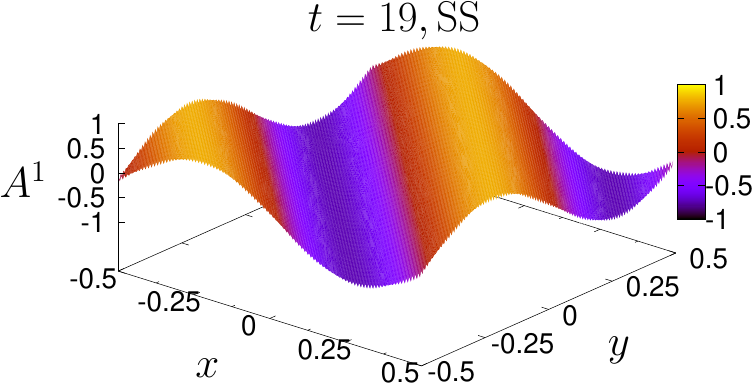}
  \end{minipage}
  \begin{minipage}{0.32\hsize}
    \includegraphics[width=\hsize]{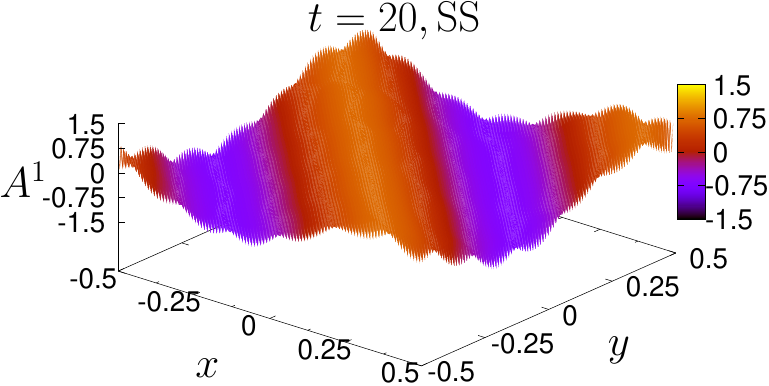}
  \end{minipage}
  \begin{minipage}{0.32\hsize}
    \includegraphics[width=\hsize]{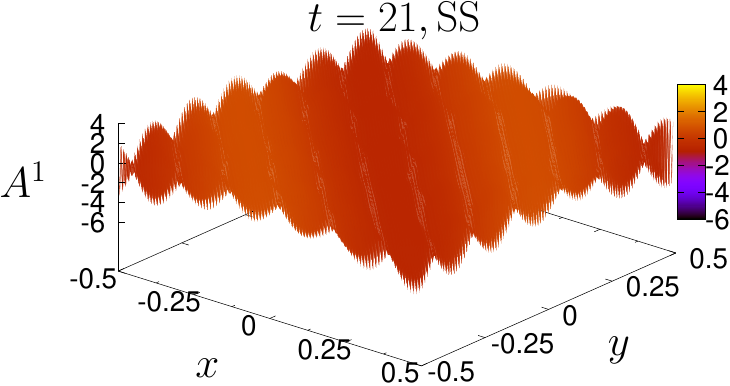}
  \end{minipage}\\
  \begin{minipage}{0.32\hsize}
    \includegraphics[width=\hsize]{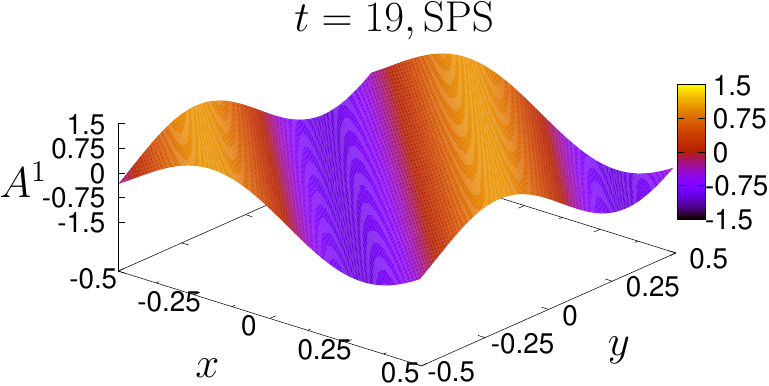}
  \end{minipage}
  \begin{minipage}{0.32\hsize}
    \includegraphics[width=\hsize]{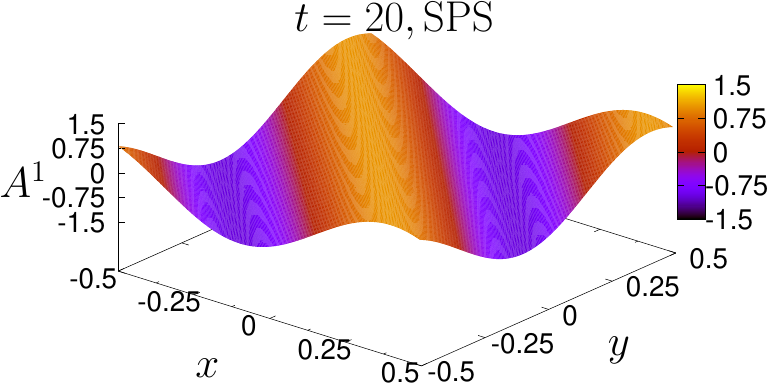}
  \end{minipage}
  \begin{minipage}{0.32\hsize}
    \includegraphics[width=\hsize]{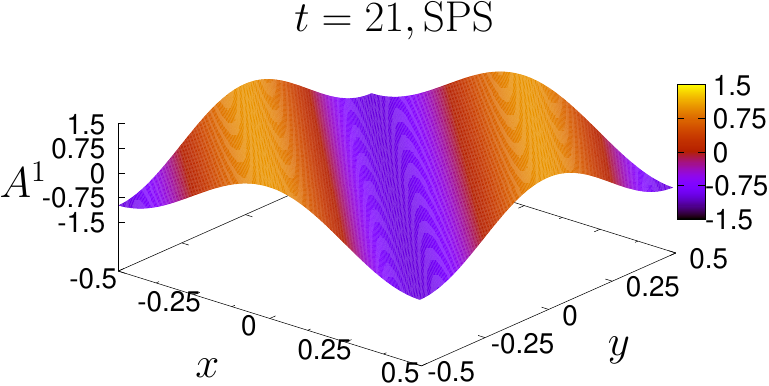}
  \end{minipage}
  \caption{
    $A^1$ with $\Delta x=\Delta y = 1/200$ at $t=19$, $20$, and $21$.
    The top panels are obtained with SS and the bottom ones with SPS.
    The left panels are at $t=19$, the middle ones at $t=20$, and the right
    ones at $t=21$.
    \label{fig:Waves_A1_t19--21}
  }
\end{figure}
Fig. \ref{fig:Waves_A1_t19--21} shows $A^1$ with $\Delta x=\Delta y=1/200$ at
$t=19$, $20$, and $21$.
At $t=20$ and $21$, we see that there are differences between the waveforms
obtained with SS and SPS.
For the waveforms obtained with SS at $t=20$ and $21$, vibrations occur.
On the other hand, at $t=19$, the behaviors of the waveforms obtained with SS
and SPS seem to be almost the same.
However, there are differences between the amplitudes obtained with SS and SPS.
\begin{figure}[htbp]
  \centering
  \begin{minipage}{0.35\hsize}
    \includegraphics[width=\hsize]{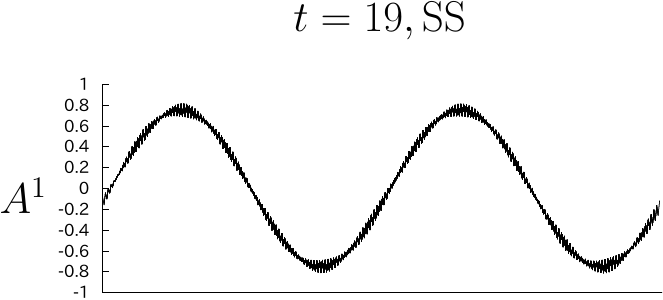}
  \end{minipage}
  \begin{minipage}{0.35\hsize}
    \includegraphics[width=\hsize]{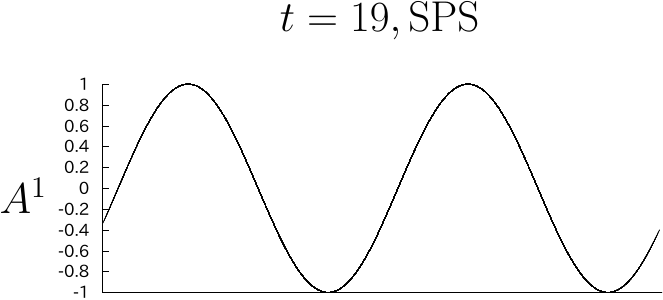}
  \end{minipage}
  \caption{
    $A^1$ with $\Delta x=\Delta y=1/200$ on $x=y$ plane.
    The left panel is for SS and the right one is for SPS.
    For SS, the sizes of $A^1$ range approximately from $-0.8$ to $0.8$, whereas
    for SPS, they range approximately from $-1$ to $1$.
    The vibration seems to occur in the waveform obtained with SS.
    \label{fig:Waves_A1_t19}
  }
\end{figure}
Fig. \ref{fig:Waves_A1_t19} shows $A^1$ with $\Delta x=\Delta y=1/200$ on the
$x=y$ plane at $t=19$.
For SS, the sizes of $A^1$ range approximately from $-0.8$ to $0.8$, whereas for
SPS, they range approximately from $-1$ to $1$.
In addition, the vibration seems to occur in the waveform obtained with SS.

\begin{figure}[htbp]
  \centering
  \begin{minipage}{0.49\hsize}
    \includegraphics[width=\hsize]{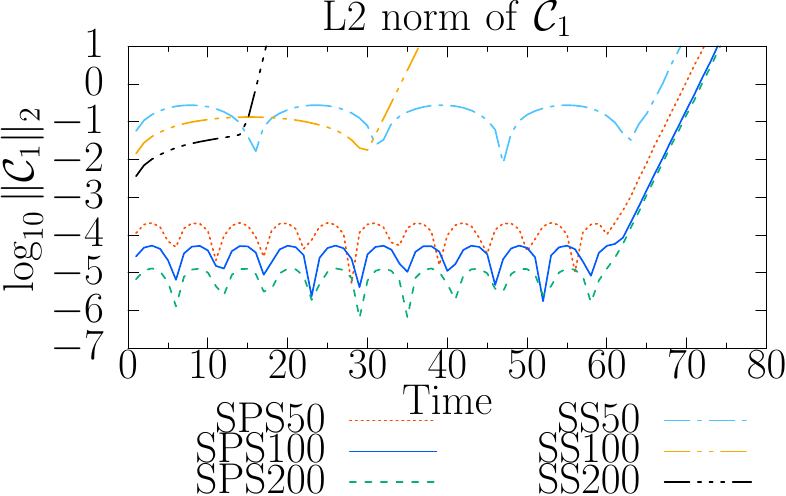}
  \end{minipage}
  \begin{minipage}{0.49\hsize}
    \includegraphics[width=\hsize]{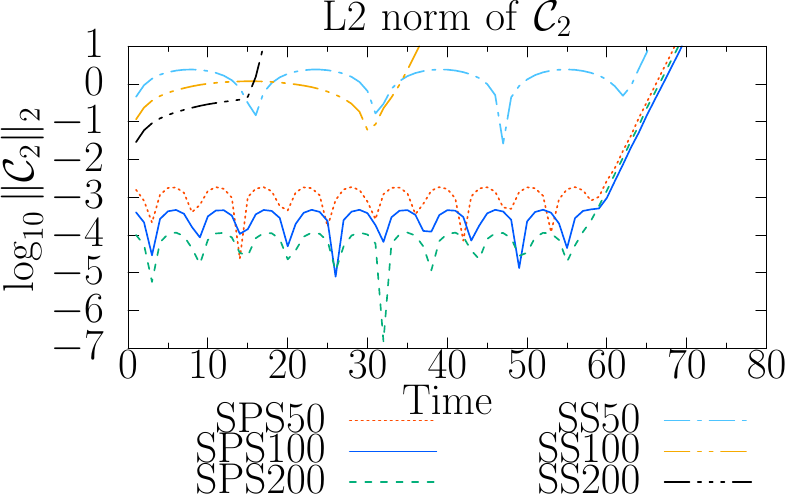}
  \end{minipage}
  \caption{
    L2 norm of $\mathcal{C}_1$ and $\mathcal{C}_2$.
    The horizontal axis indicates time and the vertical one is the
    $\log_{10}$ of the L2 norm value.
    The left panel is for $\mathcal{C}_1$ and the right one is for
    $\mathcal{C}_2$.
    The dotted line is for SPS and $\Delta x=1/50$, the solid line is for SPS
    and $\Delta x=1/100$, the dashed line is for SPS and $\Delta x=1/200$, the
    dashed dotted line is for SS and $\Delta x=1/50$, the dashed double-dotted
    line is for SS and $\Delta x=1/100$, and the dashed triple-dotted line is
    for SS and $\Delta x=1/200$.
    \label{fig:Const1-2}
  }
\end{figure}
Fig. \ref{fig:Const1-2} shows the L2 norm of $\mathcal{C}_1$ and
$\mathcal{C}_2$.
\begin{figure}[htbp]
  \centering
  \begin{minipage}{0.5\hsize}
    \includegraphics[width=\hsize]{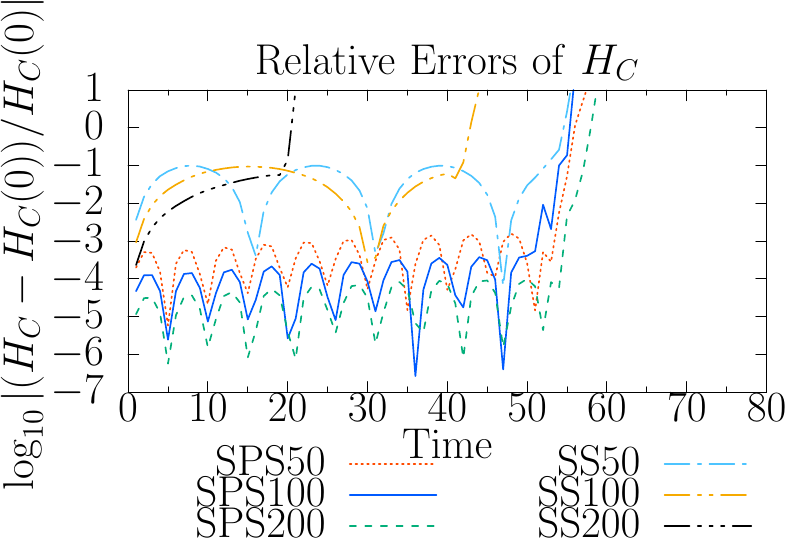}
  \end{minipage}
  \caption{
    Relative errors of $H_C$ against initial value.
    The vertical axis is $\log_{10}|(H_C-H_C(0))/H_C(0)|$.
    The others are the same as those in Fig. \ref{fig:Const1-2}.
    \label{fig:totalHC}
  }
\end{figure}
Fig. \ref{fig:totalHC} shows the relative errors of the total Hamiltonian $H_C$
against that initial value.
We see that all of the constraint values for SS are larger than those for SPS.
From the numerical results in Figs. \ref{fig:Const1-2} and \ref{fig:totalHC}, we
see that the simulations using SPS continue until about $t=55$ for all grids.
On the other hand, for SS, we see that the calculation time decreases as the
number of grids increases.
Specifically, we see that doubling the number of grids reduces the calculation
time by almost half for SS.

Setting the CFL condition to a value smaller than $1/4$ does not significantly
change the results in Figs. \ref{fig:Waves_A1_t19--21}--\ref{fig:totalHC}.

\section{Conclusion and discussion}

We introduced the canonical formulation for the Stueckelberg action of the
Proca equation, which is a wave equation for a vector field with constraints,
and proposed discrete equations with SPS.
For comparison with the conservation of the constraints at the discrete level,
we also derived other discrete equations using SS.
Numerical calculations were performed using the discrete systems with SS and
SPS, and it was shown that the results obtained with SPS are better since the
variations of the values from the initial values of the constraints
$\mathcal{C}_1$, $\mathcal{C}_2$, and $H_C$ are all smaller than those with the
obtained SS.

In Fig. \ref{fig:Const1-2}, doubling the number of grids in SS appears to reduce
the calculation time by half.
Since the CFL condition is constant, doubling the number of grids doubles the
numerical calculation time steps.
Therefore, the discretization error due to the time evolution of
$\mathcal{C}_2{}^{(\ell)}_{(\boldsymbol{k})}$ \eqref{eq:Const2Derive_SS} not
being satisfied at the discrete level has become the main discretization error.

In this study, we dealt with linear equations. A future work is to investigate
the stability and convergence of numerical calculations for nonlinear wave
equations of vector fields.

\section*{Acknowledgments}
The author was partially supported by JSPS KAKENHI Grant Numbers 24K06855
and 24K06856.


\end{document}